\documentclass[12pt]{article}
\usepackage{amsmath,amssymb,theorem}
\usepackage{graphicx}
\usepackage{color}
\usepackage{enumerate, url}
\newtheorem{thm}{Theorem}[section]
\newtheorem{conj}[thm]{Conjecture}
\newtheorem{lem}[thm]{Lemma}%[section]
\newtheorem{prop}[thm]{Proposition}%[section]
\theorembodyfont{\rmfamily}
\newtheorem{problem}[thm]{Problem}%[section]

\newcommand{\proof}{\noindent{\bf Proof.\ }}
\newcommand{\qed}{\hfill $\Box$ \bigskip}

\title{The annihilation number does not bound the $2$-domination number from the above}

\author{\small {Jun Yue, Shizhen Zhang, Yiping Zhu}\\
{\small  School of Mathematics and Statistics}\\
{\small Shandong Normal University, Jinan, Shandong 250358, P.R. China}\\
{\small Emails: yuejun06@126.com, fczsz218@126.com, yipingzhu0207@163.com}\\[2mm]
{\small Sandi Klav\v{z}ar}\\
{\small Faculty of Mathematics and Physics, University of Ljubljana, Slovenia}\\
{\small Faculty of Natural Sciences and Mathematics, University of Maribor, Slovenia}\\
{\small Institute of Mathematics, Physics and Mechanics, Ljubljana, Slovenia}\\
{\small Email: sandi.klavzar@fmf.uni-lj.si}\\[2mm]
{\small Yongtang Shi}\\
{\small Center for Combinatorics and LPMC, Nankai University, Tianjin, China}\\
{\small Email: shi@nankai.edu.cn}\\
}
\date{}
\begin{document}\maketitle
\begin{abstract}
The $2$-domination number $\gamma_2(G)$ of a graph $G$ is the minimum cardinality of a set $S\subseteq V(G)$ such that every vertex from $V(G)\setminus S$ is adjacent to at least two vertices in $S$. The annihilation number $a(G)$ is the largest integer $k$ such that the sum of the first $k$ terms of the non-decreasing degree sequence of $G$ is at most the number of its edges. It was conjectured that $\gamma_2(G) \leq a(G) +1$ holds for every connected graph $G$. The conjecture was earlier confirmed, in particular, for graphs of minimum degree $3$, for trees, and for block graphs. In this paper, we disprove the conjecture by proving that the $2$-domination number can be arbitrarily larger than the  annihilation number. On the positive side we prove the conjectured bound for a large subclass of bipartite, connected cacti, thus generalizing a result of Jakovac from [Discrete Appl.\ Math.\ 260 (2019) 178--187].
\end{abstract}

\medskip\noindent
\textbf{Keywords:} $2$-dominating set; $2$-domination number; annihilation number; cactus graph

\medskip\noindent
\textbf{AMS Math.\ Subj.\ Class.\ (2010)}: 05C69, 05C07

%%%%%%%%%%%%%%%%%%%%%%%%%%%%%%%%%%%%%%%%%%%%%%%%%%%%%%
%%%%%%%%%%%%%%%%%%%%%%%%%%%%%%%%%%%%%%%%%%%%%%%%%%%%%%
\section{Introduction}
%%%%%%%%%%%%%%%%%%%%%%%%%%%%%%%%%%%%%%%%%%%%%%%%%%%%%%
%%%%%%%%%%%%%%%%%%%%%%%%%%%%%%%%%%%%%%%%%%%%%%%%%%%%%%

Let $d_1\leq \cdots \leq d_n$ be the degree ordering of a graph $G$. The {\it annihilation number} $a(G)$ is the largest integer $k$ such that $\sum\limits_{i=1}^{k}d_i \leq m(G)$. This concept  was first defined in~\cite{Pep2004}, see also~\cite{GK1994} for an earlier, closely related concept called Havel-Hakimi process. The $2$-domination number $\gamma_2(G)$ of a graph $G$ is the minimum cardinality of a set $S\subseteq V(G)$ such that every vertex from $V(G)\setminus S$ is adjacent to at least two vertices in $S$. Now, the following conjecture relating these two concepts was posed.

\begin{conj}[\cite{D,DHRY2014}]\label{2-domination}
If $G$ is a connected graph with at least $2$ vertices, then $\gamma_2(G) \leq a(G) + 1$.
\end{conj}

If $\delta(G) \geq 3$, then Caro and Roddity~\cite[Corollary 2]{CR1990} deduced from their main result that $\gamma_2(G) \le \lfloor \frac{n(G)}{2} \rfloor$. Hence, $\gamma_2(G) \leq a(G)+1$ holds for any graph $G$ with $\delta(G) \geq 3$. Desormeaux, Henning, Rall, and Yeo~\cite{DHRY2014} followed with a confirmation of the conjecture for trees (see also~\cite{lyle-2017} for another proof of it). Moreover, they have also characterized the trees that attain the equality in the conjecture. Very recently, Jakovac~\cite{J2019} proved the conjecture for block graphs. In addition he proved:

\begin{prop}[\cite{J2019}]\label{cactus1}
If $G$ is a bipartite cactus such that every edge of $G$ belongs to a cycle, then $\gamma_2(G) \leq a(G) + 1$.
\end{prop}

In Section~\ref{sec:counter} we disprove Conjecture~\ref{2-domination} by a subclass of connected cactus graphs with minimum degree $1$. The construction further shows that the gap between the $2$-domination number and the annihilation number can be arbitrarily large. Although the conjecture is wrong, it is still interesting to find classes of graphs which satisfy the conjecture. In Section~\ref{sec:larger-class-prepare} we prove several lemmas needed in the subsequent section. Then, in Section~\ref{sec:larger-class}, we show that Conjecture~\ref{2-domination} holds for bipartite connected cacti which (i) contain no sun at an outer cycle and (ii) the degree of the exit vertex of any outer $4$-cycle is at least $4$. (A sun at a cycle is obtained from the cycle by adding a pendant vertex to each of its vertices except one.) In this way we generalize Proposition~\ref{cactus1}.  We conclude the paper with some open problems while in the rest of this section definitions and concepts needed are given.

\subsection{Preliminaries}

All graphs in this paper are undirected, finite and simple. We follow~\cite{BM2008} for graph theoretical notation and terminology not defined here.

If $G = (V(G), E(G))$ is a graph, then set $n(G)=|V(G)|$ and $m(G)=|E(G)|$. A graph $G$ is {\it nontrivial} if $n(G) \geq 2$. For $v \in V(G)$, the set of its neighbors is denoted by $N_G(v)$ and called the {\it neighborhood} of $v$, and the {\it closed neighborhood} $N_G[v]$ of $v$ is $N(v)\cup \{v\}$. The {\it degree} of a vertex $v\in V(G)$ is denoted by $d_G(v)$. For a subset $S \subseteq V(G)$, we define $\sum(S,G)=\sum_{v\in S}d_G(v)$.
In the above notation we may omit the index $G$ provided that $G$ is clear from the context. A vertex $v$ of degree $1$ is a {\it leaf} while its only neighbor is called a {\it support vertex}. If $u$ has at least two neighbors which are leaves, then $u$ is referred to as a {\it strong support vertex}. The {\it minimum} and the {\it maximum degree} among the vertices of $G$ are denoted by $\delta(G)$ and $\Delta(G)$, respectively. If $X \subseteq V(G)$, then $G-X$ denotes the graph obtained from $G$ by deleting the vertices in $X$ and all edges incident with them. Moreover, if $u_1u_2\in E(G)$ and  $v_1v_2\notin E(G)$, notations $G-u_1u_2$ and $G+v_1v_2$ will be used for the graph $(V(G),E(G)-\{u_1u_2\})$ and $(V(G),E(G) \cup \{v_1v_2\})$, respectively. These $-$ and $+$ notations will also be used for sets of edges. A connected graph is a {\it cactus} if its cycles are pairwise edge-disjoint.

A vertex $v\in V(G)$ {\it dominates} the vertices contained in $N[v]$. A set $S \subseteq V(G)$ is a {\it dominating set} if each vertex of $G$ is dominated or equivalently, if $N[S] = V(G)$, where $N[S] =\Sigma_{v\in S}N[v]$ is the {\it closed neighborhood} of $S$. The {\it domination number} $\gamma(G)$ is the minimum cardinality of a dominating set of $G$. For $k \geq 1$, a {\it $k$-dominating set} of a graph $G$ is a set $S \subseteq V(G)$ such that each vertex from $V(G)\setminus S$ is adjacent to at least $k$ vertices in $S$. There always exists at least one $k$-dominating set in $G$, since $V(G)$ is clearly a $k$-dominating set. The {\it $k$-domination number} $\gamma_k(G)$ of $G$ is the minimum cardinality of a $k$-dominating set of $G$. Thus, a $1$-dominating set is a usual dominating set and hence $\gamma_1(G) = \gamma(G)$. The notion of the $k$-dominating set was introduced by Fink and Jacobson~\cite{FJ1985}, a survey on it up to 2012 can be found in~\cite{CFH2012}. It has been further investigated afterwards,~\cite{dankelmann-2019, li-2018} are a couple of recent papers. In this paper we focus on the $2$-domination number, cf.~\cite{bujtas-2018}.

$S\subseteq V(G)$ is an {\it annihilation set} of $G$ if $\sum_{v \in S}d(v) \leq m(G)$ and is an {\it optimal annihilation set} if $|S| = a(G)$. Obviously, any optimal annihilation set of a connected graph of order at least $3$ vertices contains all leaves. Assuming that $S$ is an optimal annihilation set, we denote by $d^*(G)$ the minimum vertex degree over the set $V(G)\setminus S$.\medskip

%%%%%%%%%%%%%%%%%%%%%%%%%%%%%%%%%%%%%%%%%%%%%%%%%%%%%%
%%%%%%%%%%%%%%%%%%%%%%%%%%%%%%%%%%%%%%%%%%%%%%%%%%%%%%
\section{Counterexample to Conjecture~\ref{2-domination}}
\label{sec:counter}
%%%%%%%%%%%%%%%%%%%%%%%%%%%%%%%%%%%%%%%%%%%%%%%%%%%%%%
%%%%%%%%%%%%%%%%%%%%%%%%%%%%%%%%%%%%%%%%%%%%%%%%%%%%%%

\begin{figure}[ht!]
\begin{center}
\includegraphics[bb=191 535 393 724,scale=0.8]{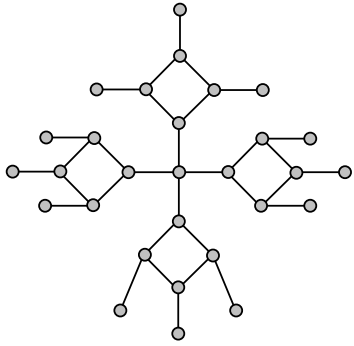}
\caption{The graph $G(4;1,1,1,1)$} \label{fig1}
\end{center}
\end{figure}

Let $t\ge 4$ and $k_1,\ldots, k_t\ge 1$. Then the graph $G(t;k_1,\ldots, k_t)$ is obtained in the following way. First, take a disjoint union of cycles $C_{3k_i+1}$, $i\in [t]$, add an additional vertex $w$, and connect $w$ with an arbitrary but fixed vertex in each of the cycles. Second, in the so far constructed graph, add a pendant vertex to each of the vertices of degree $2$. In Fig.~\ref{fig1} the graph $G(4;1,1,1,1)$ is drawn.

\begin{figure}[ht!]
\begin{center}
\includegraphics[bb=126 490 590 718,scale=0.7]{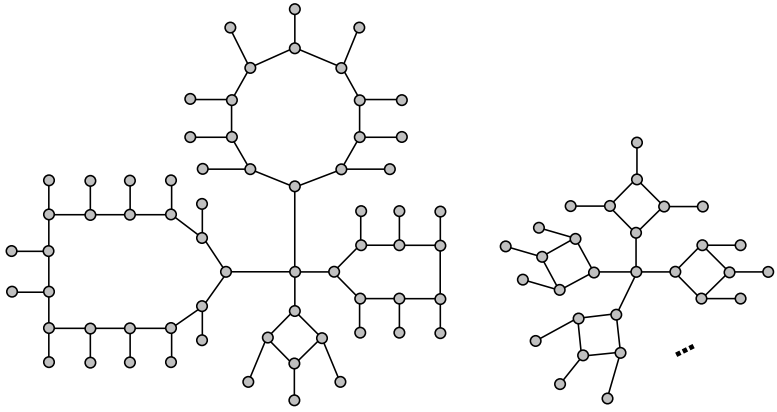}
\caption{Graphs $G(4;1,2,3,4)$ and $G(t;1,\ldots,1)$, $t\ge 4$}
\label{fig2}
\end{center}
\end{figure}

Consider first the sporadic counterexamples as shown in Fig.~\ref{fig2}. It is straightforward to verify that $a(G(4;1,2,3,4))=3+6+9+12+\left\lfloor\frac{38}{3}\right\rfloor=42$ and  $\gamma_2(G(4;1,2,3,4)) = 5+(6+3)+(9+4)+(12+5) = 44$. Hence $\gamma_2(G(4;1,2,3,4)) - a(G(4;1,2,3,4)) = 44-42 = 2>1$. Similarly, in the second example we have $a(G(t;1,\ldots,1)) = 3t + \left\lfloor\frac{8t-3t}{3}\right\rfloor
= 4t + \left\lfloor\frac{2t}{3}\right\rfloor$ and $\gamma_2(G(t;1,\ldots,1))=5t$. Therefore, $\gamma_2(G(t;1,\ldots,1)) - a(G(t;1,\ldots,1)) = 5t - (4t+\left\lfloor\frac{2t}{3}\right\rfloor)
= \left\lceil\frac{t}{3}\right\rceil$.

The above examples generalize as follows.

\begin{thm}\label{thm1.1}
Let $c_0$ be a given positive integer, $t\ge 4$, and $k_1,\ldots, k_t\ge 1$. If $t>3(c_0+1)$, then
$$\gamma_2(G(t;k_1,\ldots,k_t)) - a(G(t;k_1,\ldots,k_t)) > c_0+1\,.$$
\end{thm}

\proof
To shorten the presentation, set $G = G(t;k_1,\ldots,k_t)$ for the rest of the proof. Since to each of the constitutional cycles $C_{3k_i+1}$ of $G$ exactly $3k_i$ leaves are attached, as well as the edge to the vertex $w$ of degree $t$, we get
$$m(G) = 2\left(\sum^t_{i=1}(3k_i+1)\right)\,.$$
Hence, since each leaf of $G$ is a member of its every optimal annihilation set and all the other vertices of such a set are of degree $3$, we get
\begin{equation}
\label{eq:a-of-G}
a(G)=4\sum^t_{i=1}k_i+\left\lfloor\frac{2t}{3}\right\rfloor\,.
\end{equation}
We now claim that $\gamma_2(G)=4\sum\limits_{i=1}^t k_i+t$. Let $X$ be a $2$-dominating set of $G$ with $|X| = \gamma_2(G)$. Then every leaf of $G$ lies in $X$. Consider now a constitutional cycle $C = C_{3k_i+1}$ of $G$ and suppose that $|X\cap V(C)| \le k_i$. Then $C$ contains three consecutive vertices neither of them lying in $X$. But then the middle of these three vertices, even if being adjacent to $w$, is not $2$-dominated. If follows that
$|X\cap V(C)| \ge k_i+1$ for $i\in [t]$. Consequently, $\gamma_2(G)\geq 4\sum^t_{i=1}k_i+t$. Since on the other hand it is easy to construct a $2$-dominating set that has exactly $k_i+1$ vertices on $C$, the claim is proved. Combining the claim with~\eqref{eq:a-of-G} we conclude that
$$\gamma_2(G) - a(G) = \left(4\sum^t_{i=1}k_i+t\right) - \left(4\sum^t_{i=1}k_i+\left\lfloor\frac{2t}{3}\right\rfloor\right) = \left\lceil\frac{t}{3}\right\rceil > c_0+1\,.$$
\qed

%%%%%%%%%%%%%%%%%%%%%%%%%%%%%%%%%%%%%%%%%%%%%%%%%%%%%%
%%%%%%%%%%%%%%%%%%%%%%%%%%%%%%%%%%%%%%%%%%%%%%%%%%%%%%
\section{Some preliminary lemmas}
\label{sec:larger-class-prepare}
%%%%%%%%%%%%%%%%%%%%%%%%%%%%%%%%%%%%%%%%%%%%%%%%%%%%%%
%%%%%%%%%%%%%%%%%%%%%%%%%%%%%%%%%%%%%%%%%%%%%%%%%%%%%%

In this section, we give some lemmas to be used in the next section. They give examples of how to obtain from a graph $G$ a smaller graph $G'$, such that $\gamma_2(G') \leq a(G') + 1$ implies $\gamma_2(G) \leq a(G) + 1$. First we recall~\cite[Lemma 4]{J2019}.

\begin{lem}
\label{lem3.1}
Assume that $G$ is a graph on at least four vertices and $u\in V(G)$ a strong support vertex which is the common neighbor of pendant vertices $v_{1},\ldots,v_{\ell}\in G,~\ell \geq 2$. If $G'=G-\{u,v_{1},\ldots,v_{\ell}\}$ is a connected graph, then $\gamma_{2}(G')\leq a(G')+1$ implies $\gamma_{2}(G)\leq a(G)+1$.
\end{lem}

We proceed with new lemmas for which we define a function $f$ on a finite graph $G$ with
$$f(G) = n(G) + 3m(G) + n_1(G)\,,$$
where $n_1(G)$ denotes the number of leaves in $G$. Note that $f(G) \geq 7$ for every nontrivial, finite, connected graph $G$.

\begin{lem}\label{lem3.2}
Let $G$ be a connected graph with $n(G)\ge 3$ and which fulfils at least one of the following properties:
\begin{enumerate}
\item[(i)] $d^{*}(G)\leq 2$;
\item[(ii)] $G$ contains an induced path $P_5$ whose internal vertices are of degree $2$;
\item[(iii)] $G$ contains a pendant path $P_4$.
\end{enumerate}
Then, there exists a nontrivial connected graph $G'$ with $f(G')<f(G)$ such that $\gamma_{2}(G')\leq a(G')+1$ implies $\gamma_{2}(G)\leq a(G)+1$. Moreover, if $G$ is a connected cactus graph, then $G'$ can be chosen to be a connected cactus graph as well.
\end{lem}

\proof
If $G$ is a cycle and $e\in E(C)$, then set $G' = G-e$. If $G$ is a tree and $v$ its pendant vertex, then set $G' = G-v$. Hence in the rest of the proof we may assume that $G$ is neither a tree nor a cycle.

(i) Assume that $d^{*}(G)\leq 2$. Since $G$ is neither a tree nor a cycle, there exists a cycle $C$ in $G$ and a vertex $v \in V(C)$ with $d(v)\geq 3$. Let $e=vu \in E(C)$ and $G'=G-e$. Then $G'$ is connected, $f(G')<f(G)$ and $m(G')=m(G)-1$. The deletion of an edge does not decrease the $2$-domination number, so $\gamma_{2}(G)\leq \gamma_{2}(G')$. Consider an optimal annihilation set $S'$ of $G'$. Then $\sum(S',G')\leq m(G')=m(G)-1$. If $u,v \notin S'$, then $\sum(S',G)=\sum(S',G')\leq m(G')=m(G)-1$; if $S'$ contains exactly one of $u$ and $v$, then $\sum(S',G)=\sum(S',G')+1\leq m(G)$. In either case $a(G)\geq |S'|=a(G')$ follows. In the third case, $u,v \in S'$ and $\sum(S',G)=\sum(S',G')+2\leq m(G)+1$. Let $V_{1,2}$ denote the set of vertices which have degree $1$ or $2$ in $G$. Then $\sum(V_{1,2},G)\geq m(G)+1$ because $d^{*}(G)\leq 2$. Since $d(v)\geq 3$, we have $\sum(V_{1,2}\cup\{v\},G)\geq m(G)+4$, and then there is a vertex $v^{*}\in V_{1,2}$ which is not contained in $S'$. If $v$ is replaced with $v^{*}$ in $S'$, then we get a new annihilation set $S$ with $\sum(S,G)\leq \sum(S',G)-1\leq m(G)$. This proves $a(G)\geq |S'|=a(G')$ and then $\gamma_{2}(G)\leq\gamma_{2}(G')\leq a(G')+1\leq a(G)+1 $.

\medskip
As we have just proved the statement under the assumption (i), we can assume that $d^{*}(G)\geq 3$ in the sequel of the proof.

\medskip
(ii) Let $vu_{1}u_{2}u_{3}w$ be an induced path $P_5$ such that $d_G(u_{1}) = d_G(u_{2}) = d_G(u_{3})=2$. Set $G'=G-\{u_{1},u_{2},u_{3}\}+vw$. Observe that $n(G')=n(G)-3$, $m(G')=m(G)-3$, $n_{1}(G')=n_{1}(G)$ and hence $f(G')=f(G)-12$. Let $D'$ be a minimum $2$-dominating set of $G'$ and define $D$ as follows:
\begin{itemize}
\item[(a)] $D=D'\cup \{u_{2}\}$; if $v,w\in D'$,
\item[(b)] $D=D'\cup \{u_{1},u_{3}\}$; otherwise.
\end{itemize}
In either case, $D$ is a $2$-dominating set in $G$. Hence, $\gamma_{2}(G)\leq \gamma_{2}(G')+2$. Pick an optimal annihilation set $S'$ of $G'$. Since $d_{G}(v)=d_{G'}(v)$ and $d_{G}(w)=d_{G'}(w)$ we have $\sum(S',G)=\sum(S',G')\leq m(G')=m(G)-3$. Our assumption $d^{*}(G)\geq 3$ implies that every vertex $v$ with $d(v)\leq 2$ is contained in every optimal annihilation set of $G$. Hence, either $S'\cup\{u_{1},u_{2},u_3\}$ is an optimal annihilation set of $G$ and $a(G)\geq a(G')+2$, or there is a vertex $v^{*}\in S'$ with $d(v^{*})\geq 3$. In the latter case, consider $S=(S'-\{v^{*}\})\cup \{u_{1},u_{2},u_3\}$, and observe that $\sum(S,G)\leq \sum(S',G)-3+3\times2\leq m(G)$. Therefore, $a(G)\geq |S|\geq|S'|+2=a(G')+2$. Then $\gamma_{2}(G)\leq \gamma_{2}(G')+2 \leq a(G')+1+2 \leq a(G)+1$. This proves the statement under the assumption (ii).

\medskip
(iii) Let $u_{1}u_{2}u_{3}v$ be a pendant path $P_4$ of $G$ such that $d_G(u_{1})=1$ and $d_G(u_{2}) = d_G(u_{3}) = 2$. Since $G$ is connected, $G'=G-\{u_{1},u_{2},u_{3}\}$ is also connected, and we have $m(G')=m(G)-3$ and $f(G')<f(G)$. Let $D'$ be a minimum $2$-dominating set of $G'$. Then $D=D'\cup \{u_{1},u_{3}\}$ is $2$-dominating set of $G$. Thus, $\gamma_{2}(G)\leq \gamma_{2}(G')+2$. Next, we choose an optimal annihilation set $S'$ in $G'$. Since we have already proved (ii), we may assume that $d_G(v) \geq 3$. Consider now the following two cases. If $v \notin S'$, then $\sum(S',G)=\sum(S',G')$, and $\sum(S',G)\leq m(G')=m(G)-3$. Hence, $S=S'\cup\{u_{1},u_{2}\}$ satisfies $\sum(S,G)=\sum(S',G)+3\leq m(G)$, and $a(G)\geq |S|=|S'|+2=a(G')+2$. Then $\gamma_{2}(G)\leq\gamma_{2}(G')+2\leq a(G')+1+2\leq a(G)+1$. In the second case assume $v\in S'$. Then, $\sum(S',G)=\sum(S',G')+1\leq m(G')+1=m(G)-2$. We define $S=(S'-\{v\})\cup\{u_{1},u_{2},u_{3}\}$ and observe that $\sum(S,G)=\sum(S',G)-d(v)+5\leq m(G)-2-3+5\leq m(G)$. Hence, $S$ is an annihilation set in $G$ and we may conclude $a(G)\geq |S|=|S'|+2=a(G')+2$. So $\gamma_{2}(G)\leq\gamma_{2}(G')+2\leq a(G')+1+2\leq a(G)+1$.

To complete the proof note that all the above transformations result in a connected cactus graph $G'$, if $G$ is of the same type.
\qed

\begin{lem}\label{lem3.4}
Let $w$ be a vertex of a nontrivial, connected graph $H$ and let $v$ be a vertex of a tree $T$ with radius at least $3$, where $V(H)\cap V(T)=\emptyset$. If $G$ is the graph obtained from the $H$ and $T$ by identifying $w$ and $u$, then there exists a connected graph $G'$ with $f(G')< f(G)$ such that $\gamma_{2}(G')\leq a(G')+1$ implies $\gamma_{2}(G)\leq a(G)+1$.
\end{lem}

\proof
By Lemma~\ref{lem3.2}(i) we may suppose throughout the proof that $d^{*}(G)\geq 3$. Let $v_{1}\in V(T)$ be a vertex of $T$ at the maximum distance from $u$. Since $T$ has radius at least $3$, we have $d_T(u,v_{1}) \geq 3$. Let $v_{1}$, $v_{2}$, $v_{3}$,  $v_{4}$ be the first vertices on the shortest $v_1,u$-path in $T$ (and also in $G$). Since $d_T(u,v_{1}) \geq 3$ we infer that $v_{i} \ne u$ for  $i\in [3]$.

If $d(v_{2})\geq 3$, then $v_2$ is a strong support vertex by the assumption on $d(u,v_1)$. Then lemma holds by Lemma~\ref{lem3.1}. Hence assume in the rest  that $d(v_{2})=2$. Let $N(v_{3}) = \{w_{1}=v_2, w_2 =v_4 , \ldots, w_{t}\}$ and consider the graph $G'=G-\{v_{1},v_{2},v_{3}\} + \{w_{3}v_{4}, \ldots, w_{t}v_{4}\}$. The graph $G'$ is connected, $m(G')=m(G)-3$, and $f(G')< f(G)$. If $D'$ is a minimum $2$-dominating set of $G'$, then $D=D'\cup \{v_{1},v_{3}\}$ is $2$-dominating set of $G$. Thus, $\gamma_{2}(G)\leq |D|=|D'|+2=\gamma_{2}(G')+2$. Let $S'$ be an optimal annihilation set in $G'$. Then $\sum(S',G') \leq m(G') = m(G)-3$. Consider $S=S'\cup\{v_{1},v_{2}\}$. Then $\sum(S,G)=\sum(S',G)+d(v_{1})+d(v_{2})\leq \sum(S',G')+3\leq m(G')-3+3=m(G)$. This gives $a(G)\geq |S|=|S'|+2=a(G')+2$, and so $\gamma_{2}(G)\leq\gamma_{2}(G')+2\leq a(G')+1+2\leq a(G)+1$.
\qed

The {\it subdivided star} $S_s(K_{1,s+t})$, $s\ge 2$, $t\ge 0$, is the graph on $2s+t+1$ vertices which is constructed by subdividing $s$ edges of the star $K_{1,s+t}$ exactly once.

\begin{lem}\label{lem3.5}
Let $C$ be a cycle in a connected graph $H$ and let $w$ be a vertex of $C$ of degree $2$. If $G$ is the graph obtained from $H$ and $S_s(K_{1,s+t})$ by identifying $w$ with the center of $S_s(K_{1,s+t})$, then there exists a nontrivial connected graph $G'$ with $f(G')< f(G)$ such that $\gamma_{2}(G')\leq a(G')+1$ implies $\gamma_{2}(G)\leq a(G)+1$.
\end{lem}

\proof
Set $F = S_s(K_{1,s+t})$ and let $u$ be the center of $F$. Let $uv'_iv_i$, $i\in [s]$, be the pendant paths attached to $u$, and let $w_i$, $i\in [t]$, be the leafs adjacent to $u$, so that $V(F) = \{u, v_1,\ldots, v_s, v'_1,\ldots, v'_s, w_1, \ldots, w_t\}$. If $G'=G-V(F)$, then $G'$ is a connected cactus graph with $m(G')=m(G)-2s-t-2$ and $f(G')< f(G)$. If $D'$ is a minimum $2$-dominating set of $G'$, then $D=D' \cup \{u,w_{1},\ldots,w_{t},v_1,\ldots,v_s\}$ is a  $2$-dominating set of $G$. Thus, $\gamma_{2}(G)\leq |D|=|D'|+s+t+1= \gamma_{2}(G')+s+t+1$. Next, let $S'$ be an optimal annihilation set in $G'$. Then $\sum(S',G) \leq \sum(S',G')+2\leq m(G')+2=m(G)-2s-t$. Consider now $S=S'\cup \{v'_{1},w_{1},\ldots,w_{t},v_1,\ldots,v_s\}$. Then $\sum(S,G) \leq \sum(S',G')+s+t+4\leq m(G)-2s-t+s+t+2=m(G)-s+2 \leq m(G)$. This proves $a(G)\geq |S|=|S'|+s+t+1=a(G')+s+t+1$ and then $\gamma_{2}(G)\leq\gamma_{2}(G')+s+t+1\leq a(G')+1+s+t+1\leq a(G)+1$.
\qed

%%%%%%%%%%%%%%%%%%%%%%%%%%%%%%%%%%%%%%%%%%%%%%%%%%%%%%
%%%%%%%%%%%%%%%%%%%%%%%%%%%%%%%%%%%%%%%%%%%%%%%%%%%%%%
\section{A class of cacti for which Conjecture~\ref{2-domination} holds}
\label{sec:larger-class}
%%%%%%%%%%%%%%%%%%%%%%%%%%%%%%%%%%%%%%%%%%%%%%%%%%%%%%
%%%%%%%%%%%%%%%%%%%%%%%%%%%%%%%%%%%%%%%%%%%%%%%%%%%%%%

\begin{figure}[ht!]
\begin{center}
\includegraphics[bb=174 503 408 721,scale=0.7]{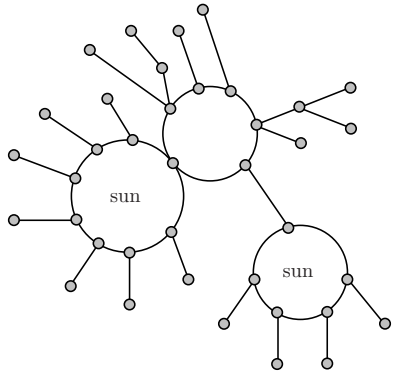}
\caption{A cactus with two suns}
\label{fig5}
\end{center}
\end{figure}

If $H_1$ and $H_2$ are subgraphs of a graph $G$, then the distance $d_G(H_1,H_2)$ between $H_1$ and $H_2$ is defined as $\min \{d_G(u,v):\ u \in V(H_1), v\in V(H_2)\}$, where $d_G(u,v)$ is the standard distance between vertices $u$ and $v$. Let $C$ and $C'$ be cycles of a cactus graph $G$. If $x \in V(C)$ and $x' \in V(C')$ such that $d_G(x,x') = d_G(C, C')$, then we say that $x$ and $x'$ are {\it exit vertices} of cycles $C$ and $C'$, respectively. A cycle of $G$ is said to be an {\it outer cycle} if it has at most one exit vertex. In the case that $G$ is unicyclic, then we also declare its cycle to be outer. Hence, if a cactus graph is not a tree, then it contains at least one outer cycle. We say that there is a {\em sun} at an outer cycle of a cactus if at all of its vertices, but at the exit vertex, there is exactly one pendant vertex attached. In Fig.~\ref{fig5} a cactus that contains two suns is drawn.

With the above terminology in hands the main result of this section reads as follows.

\begin{thm}\label{thm1.2}
Let $G$ be a connected, bipartite cactus. If $G$ contains no sun at an outer cycle, and the exit vertex of every outer $4$-cycle is of degree at least  $4$, then $\gamma_2(G) \leq a(G)+1$.
\end{thm}

\proof
We proceed by induction on the value of the function $f$ defined in the previous section. For $f(G) = 7$ we have $G \cong K_2$, and $\gamma_2(K_2) = 2 = a(K_2) + 1$. For the inductive hypothesis, let $f(G) \geq 8$ and assume that for every nontrivial graph $G'$ with $f(G') < f(G)$, we have $\gamma_2(G') \leq a(G') + 1$, where $G'$ and $G$ are connected, bipartite cactus graphs. If $G$ is a tree, then the result follows. Also, if $G$ is a cycle, then the statement is true. Thus, we may suppose that $G$ contains at least one cycle as a proper subgraph. We denote with $C_\ell$, where $\ell \geq 4$ is an even number, an outer cycle of $G$, and with $x$ the exit vertex of $C_\ell$.

In the rest of the proof we will consider subgraphs $G'$ formed from $G$ by removing a set of vertices or edges and adding edges in such a way that $f(G') < f(G)$ will hold and such that $G'$ will fulfil the assumptions of the theorem. Also, throughout the proof, $D'$ will denote a minimum $2$-domination set of $G'$, and $S'$ an optimal annihilation set in $G'$. We are going to construct a $2$-domination set $D=D' \cup D''$ and an annihilation set $S=S'\cup S''$ in $G$ that will satisfy $|D''|=|S''|=s$. Applying our inductive hypothesis to $G'$, we will estimate that $\gamma_2(G')\leq a(G')+1$ and consecutively $\gamma_2(G)\leq\gamma_2(G')+s\leq a(G')+s+1\leq a(G)+1$. In this way the theorem will be proved.

\medskip\noindent
{\bf Case 1}: All vertices from $V(C_{\ell})\setminus \{x\}$ have degree 2. \\
Let $C_\ell=x,v_1,\ldots,v_{\ell-1},x$. Set $G'=G- \{v_2,\ldots,v_{\ell-2}\}$. Then $m(G') = m(G) -(\ell-2)$. Since $d_{G'}(v_1)=d_{G'}(v_{\ell-1})=1$, both vertices $v_1$ and $v_{\ell-1}$ belong to $D'$. Set $D=D'\cup \{v_3,v_5,\ldots,v_{\ell-3}\}$.Then $D$ is a $2$-dominating set of $G$ and hence $\gamma_2(G) \leq |D| = |D'|+ \frac{\ell-4}{2} =\gamma_2(G')+\frac{\ell-4}{2}$. Since $d_{G'}(v_1)=d_{G'}(v_{\ell-1})=1$, then $v_1$ and $v_\ell$ are also both in $S'$. Set $S=S'\cup \{v_3,v_5,\ldots,v_{\ell-3}\}$, Then $\sum(S,G)\leq \sum(S',G')+2+2\frac{\ell-4}{2} \leq m(G')+2 + 2\frac{\ell-4}{2} \leq m(G)-\ell+2+(\ell-4) < m(G)$. It follows that $a(G) \geq a(G')+\frac{\ell-4}{2}$. So $ \gamma_{2}(G)\leq \gamma_{2}(G')+1 \leq a(G')+1+\frac{\ell-4}{2} \leq a(G)+1$.

\medskip\noindent
{\bf Case 2}: $V(C_{\ell})\setminus \{x\}$ contains a vertex of degree at least $3$.\\
Since $V(C_{\ell})\setminus \{x\}$ contains some vertices of degree at least $3$, and $C_{\ell}$ is an outer cycle, there are trees attached to these vertices. We root each of these trees in the vertex of the tree that lies in $V(C_{\ell})$. Amongst these trees select a tree $T$ such that $T$ has the largest height among the trees, where the height of $T$ is $\max\{d(u,v):\ u=V(C_{\ell})\cap V(T), v\in V(T)\}$. Denote the height ot $T$ with $h$, and let $u = V(T)\cap V(C_{\ell})$.

\medskip\noindent
{\bf Subcase 2.1}: $h \geq 3$.\\
Since $h \geq 3$, there exists a leaf $v \in V(T)$ such that $d(u, v) = h \geq 3$. By Lemma \ref{lem3.4} and our inductive hypothesis, the theorem holds.

\medskip\noindent
{\bf Subcase 2.2}: $h=2$. \\
We consider Cases (a), (b), (c), (d), and (e) which are schematically presented in Fig.~\ref{fig2}. All the other cases can be proved with the help of Lemma~\ref{lem3.1} and Lemma~\ref{lem3.2}(ii).

\begin{figure}[ht!]
\begin{center}
\includegraphics[bb=175 308 546 726,scale=1]{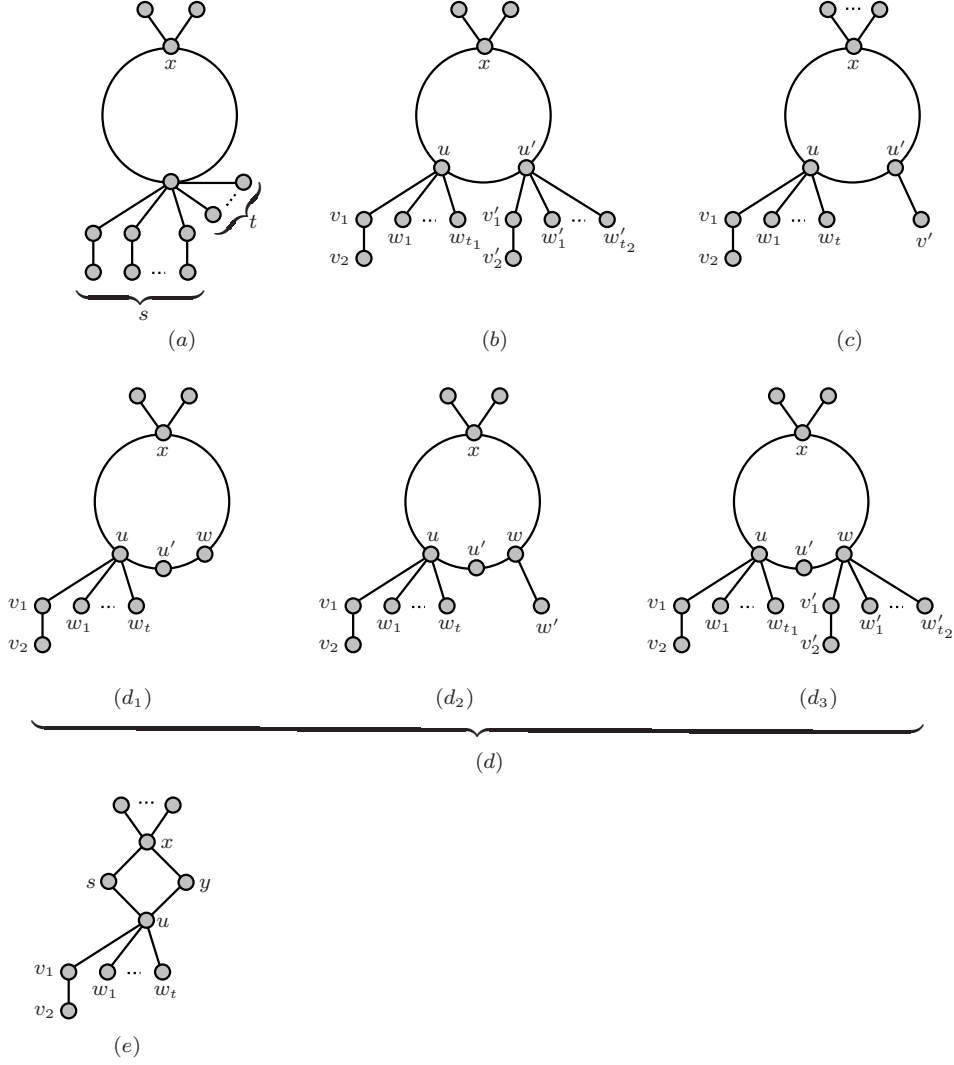}
\caption{The subcases (a), (b), (c), (d), and (e) for the case $h=2$ from the proof} \label{fig3}
\end{center}
\end{figure}

{\bf Case (a)}: In this case, we have a subdivided star $S_s(K_{1,s+t})$, ($s\geq 2 ~\text{and}~ t\geq 0$), attached to the vertex $u$ in $C_\ell$. By Lemma~\ref{lem3.5} and our inductive hypothesis for $G'=G-V(S_s(K_{1,s+t}))$, the result holds.

In the following cases we will only consider  subdivided stars with $s=1$ and $t \geq 0$, that is, the subdivided star $S_1(K_{1,1+t})$. Set $V(S_1(K_{1,1+t}))=\{u,v_1,v_2,w_1,\ldots,w_t\}$, where $u$ is the vertex of degree $t+1$, $w_1,\ldots,w_t$ are leaves adjacent to $u$, and $uv_1v_2$ is the pendant path of length $2$.

{\bf Case (b)}: In this case there are subdivided stars $S_1(K_{1,1+t_1})$ and $S_1(K_{1,1+t_2})$ with adjacent respective roots $u$ and $u'$ on $C_\ell$. Let $V(S_1(K_{1,1+t_1}))=\{u,v_1,v_2,w_1,\ldots,w_{t_1}\}$ and $V(S_1(K_{1,1+t_2}))=\{u',v'_1,v'_2,w'_1,\ldots,w'_{t_2}\}$. Set
$$G' = G- \left[ V(S_1(K_{1,1+t_1})) \cup (V(S'_1(K_{1,1+t_2})) - \{u'\})\right]\,.$$
Then $m(G')=m(G)-(t_1+t_2+6)$. Since $d_{G'}(u')=1$, we have $u' \in D'$. Set further $D=D'\cup \{u,v_2,v'_2,w_1,\ldots,w_{t_1},w'_1,\ldots,w'_{t_2}\}$. Since $D$ is a $2$-domination set of $G$ we get $\gamma_2(G)\leq|D|=|D'|+t_1+2+t_2+1=\gamma_2(G')+t_1+t_2+3$. Let $y$ be the neighbor of $u$ on $C_\ell$ different from $u'$. (Note that $y$ may be $x$.). We now consider four subcases with respect to whether $y$ and $u'$ belong to $S'$.

If $y \notin S'$ and $u'\notin S'$, we have $\sum(S',G)=\sum(S',G')\leq m(G')=m(G)-(t_1+t_2+6)$. Let $S=S'\cup\{w_1,\ldots,w_{t_1},w'_1,\ldots,w'_{t_2},v_1,v_2,v'_2\}$. Then $\sum(S,G)= \sum(S',G)+d(w_1)+\cdots+d(w_{t_1})+d(w'_1)+\cdots+d(w'_{t_2})+d(v_1)+d(v_2)+d(v'_2)\leq m(G)-(t_1+t_2+6)+t_1+t_2+2+1+1=m(G)-2<m(G)$, and we have $a(G)\geq |S'|=|S'|+t_1+t_2+3=a(G')+t_1+t_2+3$.

If $y \in S'$ and $u' \notin S'$, we have $\sum(S',G)=\sum(S',G')+1\leq m(G')+1=m(G)-(t_1+t_2+5)$. Let $S=S'\cup\{w_1,\ldots,w_{t_1},w'_1,\ldots,w'_{t_2},v_1,v_2,v'_2\}$. Then $\sum(S,G)= \sum(S',G)+d(w_1)+\cdots +d(w_{t_1})+d(w'_1)+\cdots+d(w'_{t_2})+d(v_1)+d(v_2)+d(v'_2)\leq m(G)-(t_1+t_2+5)+t_1+t_2+2+1+1=m(G)-1<m(G)$, and we have $a(G)\geq |S'|=|S'|+t_1+t_2+3=a(G')+t_1+t_2+3$.

If $y \notin S'$ and $u'\in S'$, we have $\sum(S',G)=\sum(S',G')+2+t_2\leq m(G')+2+t_2=m(G)-t_1-4$. Let $S=(S'-\{u'\})\cup\{w_1,\ldots,w_{t_1},w'_1,\ldots,w'_{t_2},v_1,v_2,v'_2,v'_1\}$. Then $\sum(S,G)= \sum(S',G)-d(u')+d(w_1)+\cdots+d(w_{t_1})+d(w'_1)+\cdots+d(w'_{t_2})+d(v_1)+d(v_2)+d(v'_2)+d(v'_1)\leq m(G)-t_1-4-t_2-3+t_1+t_2+2+1+1+2=m(G)-1<m(G)$, and we have $a(G)\geq |S|=|S'|+t_1+t_2+3=a(G')+t_1+t_2+3$.

If $y \in S'$ and $u' \in S'$, we have $\sum(S',G)=\sum(S',G')+3+t_2\leq m(G')+3+t_2=m(G)-t_1-3$. Let $S=(S'-\{u'\})\cup\{w_1,\ldots,w_{t_1},w'_1,\ldots,w'_{t_2},v_1,v_2,v'_2,v'_1\}$. Then $\sum(S,G)=\sum(S',G)-d(u')+d(w_1)+\cdots+d(w_{t_1})+d(w'_1)+\cdots+d(w'_{t_2})+d(v_1)+d(v_2)+d(v'_2)+d(v'_1)\leq m(G)-t_1-3-t_2-3+t_1+t_2+2+1+1+2=m(G)$, and we have $a(G)\geq |S'|=|S'|+t_1+t_2+3=a(G')+t_1+t_2+3$.

{\bf Case (c)}: In this case there exists a subdivided star $S_1(K_{1,1+t})$ whose vertex $u$ on $C_\ell$ has a neighbor $u'$ on $C_\ell$ with an attended pendant vertex $v'$.

If $\ell \geq 6$, then let $G'=G-(V(S_1(K_{1,1+t}))\cup \{u',v'\})$. Then $m(G')=m(G)-(t+6)$. Setting  $D=D '\cup \{u,v_2,v',w_1,\ldots,w_t\}$ we have $\gamma_2(G)\leq |D|=|D'|+t+3=\gamma_2(G')+(t+3)$. Independently of whether the neighbors of $u$ and $u'$ in $G'$ are inside $S'$ or not, we have $\Sigma(S',G)=\Sigma(S',G')+2 \leq m(G')+2=m(G)-(t+4)$. Let $S=S'\cup \{v_1,v_2,v',w_1,\ldots,w_t\}$. Then $\Sigma(S,G)= \Sigma(S',G)+(1+1+2+t)\leq m(G)-(t+4)+(t+4)=m(G)$, and we have $a(G)\geq |S|=|S'|+(t+3)=a(G')+(t+3)$.

Suppose now that $\ell = 4$ and let $C_\ell=x,u,u',y$. If $d(y) \geq 3$, then we can proceed as in the above case $\ell \geq 6$. Suppose next that $d(y)=2$. Setting $G'=G-(V(C_\ell)\cup V(S_1(K_{1,1+t})))$ we have $m(G')=m(G)-(t+7)$. Let $D=D'\cup \{u,v_2,y,v',w_1,\ldots,w_t\}$, and hence $\gamma_2(G)\leq |D|=|D'|+t+4=\gamma_2(G')+(t+4)$. If $x \notin S'$, then set $S=S'\cup \{v_1,v_2,y,v',w_1,\ldots,w_t\}$. Then $\Sigma(S,G)= \Sigma(S',G)+(2+3+1+t)\leq m(G)-(t+7)+(t+6)<m(G)$. If $x \in S'$, then set $S=(S'\setminus x)\cup \{u'\}\cup \{v_1,v_2,y,v',w_1,\ldots,w_t\}$. Then $\Sigma(S,G)= \Sigma(S',G)+2-d(x)+3+(2+3+1+t)\leq m(G)-(t+7)-1+(t+6)=m(G)$. So we have $a(G)\geq |S|=|S'|+(t+4)=a(G')+(t+4)$.

{\bf Case (d)}: In this case we have a subdivided star $S_1(K_{1,1+t})$ such that its vertex $u$ on $C_\ell$, has a neighbor $u'$ on $C_\ell$ of degree $2$. We consider thres subcases.

{\bf Case (d1)}: In this subcase $u'$ has another neighbor $w \in V(C_\ell)$ such that $d(w)=2$. If $\ell \geq 6$, select $G'=G-(V(S_1(K_{1,1+t}))\cup \{u'\})$. Then $ m(G')=m(G)-(t+5)$. Since $d_{G'}(w)$=1, we much have $w \in D'$. Let $D=D'\cup \{u,v_2,w_1,\ldots,w_t\}$. Since $D$ is a $2$-dominating set of $G$, we get $\gamma_2(G)\leq |D|=|D'|+(t+2)=\gamma_t(G')+(t+2)$. Independently of whether the neighbors of $u$ and $u'$ in $G'$ are inside $S'$ or not, we have $\Sigma(S',G)=\Sigma(S',G')+2\leq m(G')+2=m(G)-(t+3)$. Let $S=S'\cup \{v_1,v_2,w_1,\ldots,w_t\}$. Then $\Sigma(S,G)=\Sigma(S',G)+(t+3)\leq m(G)-(t+3)+(t+3)\leq m(G)$, and we have $a(G)\geq |S|=|S'|+(t+2)=a(G')+(t+2)$.

Consider now the case that $\ell =4$. Setting $G'=G-(V(S_1(K_{1,1+t}))\cup \{u',w\})$ we have $m(G')=m(G)-(t+6)$. Let $D=D '\cup \{u,v_2,w,w_1,\ldots,w_t\}$, and hence $\gamma_2(G)\leq |D|=|D'|+t+3=\gamma_2(G')+(t+3)$. If $x \notin S'$, then set $S=S'\cup \{v_1,v_2,w,w_1,\ldots,w_t\}$. Then $\Sigma(S,G)= \Sigma(S',G)+(2+2+1+t)\leq m(G)-(t+6)+(t+5)<m(G)$. If $x \in S'$, then set $S=(S'\setminus x)\cup \{u'\}\cup \{v_1,v_2,w,w_1,\ldots,w_t\}$. Then $\Sigma(S,G)= \Sigma(S',G)+2-d(x)+2+(2+2+1+t)\leq m(G)-(t+6)+(t+5)<m(G)$. So we have $a(G)\geq |S|=|S'|+(t+3)=a(G')+(t+3)$.

{\bf Case (d2)}:
Suppose new that the other neighbor $w$ of $u'$ has a pendant vertex $w'$. If $\ell \geq 6$, then let $G'=G-(V(S_1(K_{1,1+t}))\cup \{u',w'\})$. Then $ m(G')=m(G)-(t+6)$. Since $d_{G'}(w)$=1, we have $w \in D'$. Set $D=D'\cup \{u,v_2,w',w_1,\ldots,w_t\}$. Then $D$ is a $2$-dominating set of $G$ and therefore $\gamma_2(G)\leq |D|=|D'|+(t+3)=\gamma_t(G')+(t+3)$. If $w \notin S'$, then we have $\Sigma(S',G)\leq \Sigma(S',G')+1 \leq m(G)-(t+5)$. Set next $S=S'\cup \{v_1,v_2,w',w_1,\ldots,w_t\}$. Then $\Sigma(S,G)=\Sigma(S',G)+(1+2+1+t)\leq m(G)-(t+5)+(t+4)<m(G)$, and we have $a(G)\geq |S|=|S'|+(t+3)=a(G')+(t+3)$. If $w \in S'$, we have $\Sigma(S',G)\leq \Sigma(S',G')+3 \leq m(G')+3=m(G)-3$. Setting $S=(S'-\{w\})\cup \{u',v_1,v_2,w',w_1,\ldots,w_t\}$, we have $\Sigma(S,G)=\Sigma(S',G)-3+(2+2+1+1+t)\leq m(G)-(t+3)+(t+3)=m(G)$, and we have $a(G)\geq |S|=|S'|+(t+3)=a(G')+(t+3)$.

Let now $\ell =4$. Setting $G'=G-(V(S_1(K_{1,1+t}) \cup \{u',w,w'\})$ we have  $m(G')=m(G)-(t+7)$. Let $D=D '\cup \{u,v_2,w,w',w_1,\ldots,w_t\}$, and hence $\gamma_2(G)\leq |D|=|D'|+t+4=\gamma_2(G')+(t+4)$. If $x \notin S'$, then set $S=S' \cup \{u,v_2,u',w',w_1,\ldots,w_t\}$ and $\Sigma(S,G)= \Sigma(S',G)+(2+2+1+1+t)\leq m(G)-(t+7)+(t+6)<m(G)$. If $x \in S'$, then set $S=(S'\setminus x)\cup \{w\}\cup \{u,v_2,u',w',w_1,\ldots,w_t\}$ and therefore $\Sigma(S,G)= \Sigma(S',G)+2-d(x)+3+(2+2+1+1+t)\leq m(G)-(t+7)+1+(t+6)=m(G)$. So we have $a(G)\geq |S|=|S'|+(t+4)=a(G')+(t+4)$.

{\bf Case (d3)}:
Suppose now that at the other neighbor $w$ of $u'$ we have another subdivided star $S_1(K_{1,1+t_2})$. %If $\ell \geq 6$, then
Set $G'=G-(V(S_1(K_{1,1+t_1}) \cup V(S_1(K_{1,1+t_2}) \cup \{u'\})$ and hence $ m(G')=m(G)-(t_1+t_2+8)$. Setting $D=D'\cup \{u,w,v_2,v'_2,w_1,\ldots,w_{t_1}, w'_1,\ldots,w'_{t_2}\}$ we get $\gamma_2(G)\leq|D|=|D'|+(t_1+t_2+4)=\gamma_2(G')+(t_1+t_2+4)$. Independently of whether the neighbors of $u$ and $w$ in $G'$ are inside $S'$ or not, we have $\Sigma(S',G)\leq \Sigma(S',G')+2 \leq m(G)-(t_1+t_2+6)$. Set further $S=S'\cup \{v_1,v_2,v'_1,v'_2,w_1,\ldots,w_{t_1},w'_1,\ldots,w'_{t_2}\}$, so that $\Sigma(S,G)=\Sigma(S',G)+(1+2+1+2+t_1+t_2)\leq m(G)-(t_1+t_2+6)+(t_1+t_2+6) = m(G)$, and $a(G)\geq |S|=|S'|+(t_1+t_2+4)=a(G')+(t_1+t_2+4)$.

%Let now $\ell =4$. Setting $G'=G-(V(S_1(K_{1,1+t_1}) \cup V(S_1(K_{1,1+t_2}) \cup \{u'\})$ we have  $m(G')=m(G)-(t_1+t_2+8)$. Let $D=D '\cup \{u,v_2,w,v'_2,w_1,\ldots,w_t,w'_1,\ldots,w'_{t_2}\}$, and hence $\gamma_2(G)\leq |D|=|D'|+t_1+t_2+4=\gamma_2(G')+(t_1+t_2+4)$. Independently of whether $x$ is insider $S'$ or not, set $S=S' \cup \{v_1,v_2,v'_1,v'_2,w_1,\ldots,w_t,w'_1,\ldots,w'_{t_2}\}$ and then $\Sigma(S,G)= \Sigma(S',G)+2+(2+2+1+1+t_1+t_2)\leq m(G)-(t_1+t_2+8)+2+(t_1+t_2+6)=m(G)$. So we have $a(G)\geq |S|=|S'|+(t_1+t_2+4)=a(G')+(t+4)$.

{\bf Case (e)}: Let $\ell = 4$ and let $C_4 =x,s,u,y,x$, where $u$ is attended with a subdivided star $S(K_{1,1,t})$. Setting $G'=G-V(S_1(K_{1,1+t}))$ we have $ m(G')=m(G)-(t+4)$. Since $d_G'(s)=d_G'(y)=1$, the vertices $s$ and $y$ belong to $D'$. Let $D=\{D'\setminus \{s,y\}\}\cup \{x,u,v_2,w_1,\ldots,w_{t}\}$. Then $\gamma_2(G)\leq|D|=|D'|+(t+1)=\gamma_2(G')+(t+1)$. Set further  $S=S'\cup \{v_2,w_1,\ldots,w_{t}\}$. Then $\Sigma(S,G)=\Sigma(S',G)+2+(1+t)\leq m(G)-(t+4)+(t+3) < m(G)$, and $a(G)\geq |S|=|S'|+(t+1)=a(G')+(t+1)$.

\begin{figure}[ht!]
\begin{center}
\includegraphics[bb=243 560 357 717,scale=0.8]{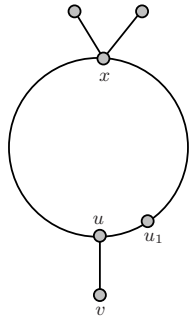}
\caption{The situation to be considered in the subcase $h=1$} \label{fig4}
\end{center}
\end{figure}

\medskip\noindent
{\bf Subcase 2.3}: $h=1$.\\
We need to consider only one case which is shown in Fig.~\ref{fig4}, because, as we have already seen in Case 2.2, all the other cases for $h = 1$ can be proved with the help of Lemma~\ref{lem3.1}.

Assume that $d(u)=3$, and at least one of its neighbors in $V(C_\ell)\setminus \{x\}$, is degree of $2$, denote it with $u_1$. Denote the child of $u$ with $v$. Setting $G'=G-\{u,v\}$ we have  $m(G')=m(G)-3$. Then $u_1\in D'$, since it is leaf in $G'$. Setting $D=D'\cup\{v\}$ we get $\gamma_{2}(G)\leq|D|=|D'|+1=\gamma_{2}(G')+1$. Independently of whether the neighbors of $u$ in $G'$ are inside $S'$ or not, we have $\sum(S',G)\leq\sum(S',G')+2\leq m(G')+2=m(G)-1$. Let $S=S'\cup\{v\}$. Then $\sum(S,G)= \sum(S',G)+d(v)\leq m(G)-1+1=m(G)$, and we have $a(G)\geq |S'|=|S'|+1=a(G')+1$.
\qed

%%%%%%%%%%%%%%%%%%%%%%%%%%%%%%%%%%%%%%%%%%%%%%%%%%%%%%
%%%%%%%%%%%%%%%%%%%%%%%%%%%%%%%%%%%%%%%%%%%%%%%%%%%%%%
\section{Concluding remarks}
%%%%%%%%%%%%%%%%%%%%%%%%%%%%%%%%%%%%%%%%%%%%%%%%%%%%%%
%%%%%%%%%%%%%%%%%%%%%%%%%%%%%%%%%%%%%%%%%%%%%%%%%%%%%%

Based on the results of this paper, the following problem is very natural.

\begin{problem}
Characterize the cactus graphs for which Conjecture~\ref{2-domination} holds true.
\end{problem}

Note that the class of cacti in question does not contain bipartite cacti as a subclass since some of the counterexamples from Section~\ref{sec:counter} are bipartite. More generally, it would be interesting to know the answer to the following:

\begin{problem}
Characterize the graphs for which Conjecture~\ref{2-domination} holds true.
\end{problem}

As we already mentioned, in~\cite{DHRY2014} trees were characterized for which the equality in Conjecture~\ref{2-domination} holds. Hence we pose:

\begin{problem}
Characterize the cactus graphs for which the equality in Conjecture~\ref{2-domination} holds. More generally, characterize the graphs with the same property.
\end{problem}

Let $\gamma_t(G)$ denote the total domination number of a graph $G$. (For an extensive information on $\gamma_t$ see the book~\cite{henning-2013}.) In~\cite{D, DHRY2014} a parallel conjecture to Conjecture~\ref{2-domination} was posed for the total domination number, that is, it was conjectured that
\begin{equation}
\label{eq:conj}
\gamma_t(G) \le a(G) + 1
\end{equation}
holds for every nontrivial connected graph $G$. This conjecture holds for graphs of minimum degree at least $3$, and has been verified for trees~\cite{desormeaux-2013} and for cactus graphs and block graphs~\cite{bujtas-2019}. The counterexamples to Conjecture~\ref{2-domination} presented in this paper are far from being counterexamples for~\eqref{eq:conj} since their total domination number is significantly smaller and, after all, the counterexamples to Conjecture~\ref{2-domination} are cactus graphs for which~\eqref{eq:conj} holds. Hence we are inclined to believe that~\eqref{eq:conj} holds true.

%%%%%%%%%%%%%%%%%%%%%%%%%%%%%%%%%%%%%%%%%%%%%%%%%%%%%%
%%%%%%%%%%%%%%%%%%%%%%%%%%%%%%%%%%%%%%%%%%%%%%%%%%%%%%
\section*{Acknowledgments}
%%%%%%%%%%%%%%%%%%%%%%%%%%%%%%%%%%%%%%%%%%%%%%%%%%%%%%
%%%%%%%%%%%%%%%%%%%%%%%%%%%%%%%%%%%%%%%%%%%%%%%%%%%%%%

Jun Yue was partially supported by the National Natural Science Foundation of China (No.\ 11626148 and 11701342) and the Natural Science Foundation of Shandong Province (No.\ ZR2016AQ01). Yongtang Shi was partially supported by the National Natural Science Foundation of China, Natural Science Foundation of Tianjin
(No.\ 17JCQNJC00300), the China-Slovenia bilateral project ``Some topics in modern graph theory" (No.\ 12-6),  Open Project Foundation of Intelligent Information Processing Key Laboratory of Shanxi
Province (No.\ CICIP2018005), and the Fundamental Research Funds for the Central Universities, Nankai University (63191516). Sandi Klav\v{z}ar acknowledges the financial support from the Slovenian Research Agency (research core funding P1-0297 and projects J1-9109, N1-0095, N1-0108).

\end{document}